# An Accurate and Efficient Algorithm for The Time-fractional Molecular Beam Epitaxy Model with Slope Selection

Lizhen Chen *, Jia Zhao †, Waixiang Cao ‡, Hong Wang §and Jiwei Zhang ¶

January 18, 2018

**Abstract**

In this paper, we propose a time-fractional molecular beam epitaxy (MBE) model with slope selection and its efficient, accurate, full discrete, linear numerical approximation. The numerical scheme utilizes the fast algorithm for the Caputo fractional derivative operator in time discretization and Fourier spectral method in spatial discretization. Refinement tests are conducted to verify the $2-\alpha$ order of time convergence, with $\alpha \in (0,1]$ the fractional order of derivative. Several numerical simulations are presented to demonstrate the accuracy and efficiency of our newly proposed scheme. By exploring the fast algorithm calculating the Caputo fractional derivative, our numerical scheme makes it practice for long time simulation of MBE coarsening, which is essential for MBE model in practice. With the proposed fractional MBE model, we observe that the scaling law for the energy decays as $O(t^{-\frac{\alpha}{3}})$ and the roughness increases as $O(t^{\frac{\alpha}{3}})$, during the coarsening dynamics with random initial condition. That is to say, the coarsening rate of MBE model could be manipulated by the fractional order $\alpha$, and it is linearly proportional to $\alpha$. This is the first time in literature to report/discover such scaling correlation. It provides a potential application field for fractional differential equations. Besides, the numerical approximation strategy proposed in this paper can be readily applied to study many classes of time-fractional and high dimensional phase field models.

## 1 Introduction

Molecular beam epitaxy (MBE) method is a broadly used approach for thin-film deposition of single crystal. And this strategy is widely applied in semiconductor manufacture.

*Beijing Computational Science Research Center, Beijing, China; Email: lzchen@csrc.ac.cn.
†Department of Mathematics & Statistics, Utah State University, Logan, UT, USA; Email: jia.zhao@usu.edu.
‡School of Mathematical Science, Beijing Normal University. E-mail: caowx@bnu.edu.cn.
§Department of Mathematics, University of South Carolina. E-mail: hwang@math.sc.edu.
¶Beijing Computational Science Research Center, Beijing, China; Email: jwzhang@csrc.ac.cn.

In recent years, MBE becomes a very important and challenging research topic in material science. In the meanwhile, many mathematical models have been developed to study the epitaxy dynamics, ranging from molecular dynamics simulations to continuum models [5, 9, 12, 13, 18, 20, 23, 27, 34, 38].

In the continuum model approach, one popular continuum model for molecular beam epitaxy is derived via energy variational approach, satisfying an energy dissipation law (i.e., thermodynamically consistent) [27]. Consider a smooth domain $\Omega$, and use $\phi(\mathbf{x}, t) : \Omega \to \mathbb{R}$ to denote the height function of MBE, and the effective free energy is given as

$$E(\phi) = \int_\Omega \Big[\frac{\varepsilon^2}{2}|\Delta \phi|^2 + f(\nabla \phi)\Big] d\Omega. \tag{1.1}$$

Here the first term represents the isotropic surface diffusion effect with $\varepsilon$ a constant controlling the surface diffusion strength, and the second term approximates the Enrlich-Schwoebel effect that the adatoms stick to the boundary from an upper terrace, contributing to the steepening of mounds in the film [7]. The evolution equation for $\phi$ could be derived via a $L^2$ gradient flow associated with the effective free energy functional $E(\phi)$, i.e. the equation for $\phi$ reads as

$$\partial_t \phi = -M \frac{\delta E}{\delta \phi}, \tag{1.2}$$

where $M$ is the mobility parameter (with $\frac{1}{M}$ proportional to the relaxation time). For simplicity of notations, we consider periodic boundary condition. If we choose $f(\nabla \phi) = -\frac{1}{2}\ln(1 + |\nabla \phi|^2)$, the corresponding equation would reduce to

$$\begin{cases} \partial_t \phi(\mathbf{x}, t) = -M\Big(\varepsilon^2 \Delta^2 \phi + \nabla \cdot \big(\frac{\nabla \phi}{1 + |\nabla \phi|^2}\big)\Big), & (\mathbf{x}, t) \in \Omega \times (0, T] \\ \phi(\mathbf{x}, 0) = \phi_0(\mathbf{x}), & \mathbf{x} \in \Omega. \end{cases} \tag{1.3}$$

and its energy dissipation rate could be calculated as

$$\frac{dE}{dt} = -\int_\Omega M\Big(\varepsilon^2 \Delta^2 \phi + \nabla \cdot \big(\frac{\nabla \phi}{1 + |\nabla \phi|^2}\big)\Big)^2 d\Omega. \tag{1.4}$$

If we choose $f(\nabla \phi) = \frac{1}{4}(|\nabla \phi|^2 - 1)^2$, the corresponding equation becomes

$$\begin{cases} \partial_t \phi(\mathbf{x}, t) = -M\Big(\varepsilon^2 \Delta^2 \phi + \nabla \cdot ((1 - |\nabla \phi|^2)\nabla \phi)\Big), & (\mathbf{x}, t) \in \Omega \times (0, T], \\ \phi(\mathbf{x}, 0) = \phi_0(\mathbf{x}), & \mathbf{x} \in \Omega. \end{cases} \tag{1.5}$$

and its energy dissipation rate can be calculated as

$$\frac{dE}{dt} = -\int_\Omega M\Big(\varepsilon^2 \Delta^2 \phi + \nabla \cdot ((1 - |\nabla \phi|^2)\nabla \phi)\Big)^2 d\Omega. \tag{1.6}$$

In addition, both models (1.3) and (1.5) conserve the total mass, i.e.

$$\frac{d}{dt}\int_\Omega \phi(\mathbf{x},t)d\Omega = 0. \tag{1.7}$$

There exists a huge amount of work in literature on investigating this MBE model (1.5) analytically and numerically [6–8, 10, 11, 17, 19, 21–23, 29–32, 35, 36, 39, 43, 45]. In the MBE model (1.5), it has been shown that the roughness (the standard deviation of the height profile) scales as $O(t^{\frac{1}{3}})$ and the effective energy scales as $O(t^{-\frac{1}{3}})$ [19, 22, 23, 39, 45]. One might wonder how the scaling rate could be manipulated. A straight forward idea is to introduce the time-fractional derivative, which takes the long-time memory behavior into the consideration.

In recent years, researchers are conducting extensive research on generalizing the integer-order partial differential equations (PDEs) into fractional PDEs to better model the anomalously diffusive effects. In particular, the computation of fractional PDE models become practical, with the development of advanced fast algorithms both in time PDEs [16, 24, 44] and space fractional differential equations [40, 41]. The understanding of fractional PDE models' behavior and their applications are unprecedented [3, 4, 33, 37]. Among the field of fractional PDEs, one specific research topic is the fractional gradient flow problems. The gradient flow problems are derived via an energetic variation, such that the derived model would satisfy an energy dissipation law, i.e., thermodynamically consistent. One essential question is whether the fractional generalization would still possess the energy dissipation property. Several pieces of work have been published to show that the spatial-fractional generalization do possess the energy dissipation property. For instance, Ainsworth and Mao show that a fractional Cahn-Hilliard equation by considering a gradient flow in the negative order Sobolev space $H^{-\beta}, \beta \in [0,1]$ obey energy dissipation law [1, 2]. Similar results are also observed in [42]. Unfortunately, it is still an open question for the time-fractional generalization, though several numerical results for time-fractional gradient flow models do indicate that they have energy dissipation properties [25, 26].

In this paper, we focus on a generalized version of the MBE model (1.5). Mainly, we consider the time-fractional MBE model with slope selection, which reads

$$\begin{cases} {}_0^C D_t^\alpha \phi(\mathbf{x},t) = -M\Big(\varepsilon^2 \Delta^2 \phi - \nabla\cdot(((|\nabla\phi|^2-1)\nabla\phi))\Big), & (\mathbf{x},t) \in \Omega\times(0,T] \\ \phi(\mathbf{x},0) = \phi_0(\mathbf{x}), \quad \mathbf{x}\in\Omega, \end{cases} \tag{1.8}$$

where $\alpha \in (0,1]$ and ${}_0^C D_t^\alpha$ is the Caputo fractional derivative of order $\alpha$ defined by [28]

$$_0^C D_t^\alpha \phi(t) := \begin{cases} \dfrac{1}{\Gamma(1-\alpha)}\displaystyle\int_0^t \dfrac{\phi_\eta(\eta)}{(t-\eta)^\alpha}d\eta, & 0 \leq \alpha \leq 1, \\ \dfrac{\partial \phi}{\partial t}, & \alpha = 1, \end{cases} \tag{1.9}$$

with $\Gamma(\cdot)$ the Gamma function. The Caputo type time-fractional differential operator provides a means to model the sub-diffusive or long time memory behavior. The main goal of this paper is to develop fast and accurate numerical schemes for the time-fractional MBE model, and to investigate its dynamic differences compared with the integer-order MBE model. Furthermore, through numerical tests, we verify that it obeys an energy dissipation law, and we investigate the energy decay rate and roughness growth rate by the coarse dynamics simulation with random initial profile. We discover that for the fractional MBE model (1.8) the roughness scales like $t^{\frac{\alpha}{3}}$ and the effective energy scales like $t^{-\frac{\alpha}{3}}$, which will be explained in details in this paper.

The outline of the paper is listed as follows. In Section 2, we present the time-fractional MBE model. In Section 3, we develop the numerical scheme. The fast algorithm for the Caputo fractional derivative operator and fourier spectral method were chosen for the time and space discretization respectively. In Section 4, we present some numerical simulations to demonstrate the accuracy and efficiency of the proposed schemes. Finally, some concluding remarks are presented in Section 5.

## 2  The time-fractional Molecular Beam Epitaxy Model

For simplicity of notations, we assume periodic boundary conditions in the rest of this paper. We consider the generalized MBE model, namely the time fractional MBE model with slope selection as follows:

$$\begin{cases} {}_{0}^{C}D_t^{\alpha}\phi(\mathbf{x},t) = -M\Big(\varepsilon^2\Delta^2\phi - \nabla\cdot\big((|\nabla\phi|^2-1)\nabla\phi\big)\Big), & (\mathbf{x},t)\in\Omega\times(0,T] \\ \phi(\mathbf{x},0) = \phi_0(\mathbf{x}), & \mathbf{x}\in\Omega, \end{cases} \quad (2.1)$$

where $\alpha \in (0,1]$ and ${}_{0}^{C}D_t^{\alpha}$ is the classical Caputo fractional derivative of order $\alpha$ defined in (1.9).

Denoting the Riemann-Liouville fractional derivative ${}_{0}^{R}D_t^{\alpha}$ as following

$$ {}_{0}^{R}D_t^{\alpha}\phi(\mathbf{x},t) = \frac{1}{\Gamma(1-\alpha)}\frac{d}{dt}\int_0^t \frac{\phi(\eta)}{(t-\eta)^{\alpha}}d\eta \quad (2.2)$$

Taking the derivative ${}_{0}^{R}D_t^{1-\alpha}$ on both side of (2.1), then we have

$$\begin{cases} D_t^1\phi(\mathbf{x},t) = -M{}_{0}^{R}D_t^{1-\alpha}\Big(\varepsilon^2\Delta^2\phi - \nabla\cdot\big((|\nabla\phi|^2-1)\nabla\phi\big)\Big), & (\mathbf{x},t)\in\Omega\times(0,T], \\ \phi(\mathbf{x},0) = \phi_0(\mathbf{x}), & \mathbf{x}\in\Omega. \end{cases} \quad (2.3)$$

The time-fractional MBE has the mass conservation property. Actually, it is easy to

calculate that

$$\frac{d}{dt}\int_\Omega \phi(\mathbf{x},t)d\Omega = -M_0^R D_t^{1-\alpha}\int_\Omega \Big(\varepsilon^2\Delta^2\phi - \nabla\cdot((|\nabla\phi|^2-1)\nabla\phi)\Big)d\Omega = 0. \qquad (2.4)$$

Unfortunately, there is no clue as whether it obeys an energy dissipation law or not. As a matter of fact, the effective energy changing rate in time could be calculated as

$$\begin{aligned}
\frac{dE}{dt} &= \int_\Omega \frac{\delta E}{\delta \phi}\frac{\delta \phi}{\delta t}d\Omega \\
&= -\int_\Omega M\Big(\varepsilon^2\Delta^2\phi - \nabla\cdot((|\nabla\phi|^2-1)\nabla\phi)\Big)\,_0^R D_t^{1-\alpha}\Big(\varepsilon^2\Delta^2\phi - \nabla\cdot((|\nabla\phi|^2-1)\nabla\phi)\Big)d\Omega.
\end{aligned} \qquad (2.5)$$

However, it is still an open question that whether the right hand side of (2.5) is non-positive or not. Thus, we resort to the numerical simulations to verify the energy law of the time-fractional MBE model with slope selection (2.1).

In the rest of this paper, we would propose accurate and efficient numerical schemes for solving the time-fractional MBE model with slope selection (2.1), and investigate its coarsening dynamics, scaling law correlations with the fractional derivative order $\alpha$.

## 3  Numerical Schemes

In this section, we propose a fully discrete spectral method for the time-fractional MBE model (2.1). Firstly, we introduce the discretization of the Caputo fractional derivative operator, and the semi-discrete scheme for the fractional MBE model. Then we use Fourier-spectral method to handle spatial discretization.

### 3.1  Time discrete scheme

Define $\Delta t = T/K$ as the grid step in time, then the time-fractional term at $t^{k+1}$ can be approximated by the following scheme [14, 15]

$$\begin{aligned}
&{}_0^C D_t^\alpha \phi(\mathbf{x},t_{k+1}) \\
&= \frac{1}{\Gamma(1-\alpha)}\sum_{j=0}^k \int_{t_j}^{t_{j+1}} \frac{\phi_s(\mathbf{x},s)}{(t_{j+1}-s)^\alpha}ds, \\
&= \frac{1}{\Gamma(1-\alpha)}\sum_{j=0}^k \frac{\phi(\mathbf{x},t_{j+1})-\phi(\mathbf{x},t_j)}{\Delta t}\int_{j\Delta t}^{(j+1)\Delta t}\int_{t_j}^{t_{j+1}}\frac{ds}{(t_{j+1}-s)^\alpha} + O(\Delta t^{2-\alpha}) \\
&= \frac{1}{\Gamma(2-\alpha)}\sum_{j=0}^k \frac{\phi(\mathbf{x},t_{k+1-j})-\phi(\mathbf{x},t_{k-j})}{\Delta t^\alpha}[(j+1)^{1-\alpha}-j^{1-\alpha}] + O(\Delta t^{2-\alpha}).
\end{aligned} \qquad (3.1)$$

It is obvious that this time discrete scheme is nonlocal, such that it is not practical for long time simulations. Especially in two or three dimensions, a scheme using this discretization will eat up the available machine memories for simulations in long time period. To overcome this difficulty, a fast algorithm for the numerical approximation of the Caputo fractional derivative ${}_0^C D_t^\alpha \phi(\mathbf{x}, t)$ of order $0 \leq \alpha \leq 1$ was proposed by [16, 44]. In this paper, we utilize this fast algorithm for evaluation of ${}_0^C D_t^\alpha \phi(\mathbf{x}, t)$ as all the numerical simulations for MBE are to investigate two dimensional and long-time dynamics. The main idea of the fast algorithm will be described in the next subsection.

Consider the fast evaluation of the Caputo fractional derivative which was proposed by [16, 44]. Before introduce the scheme, we approximate the kernel $t^{-\beta}(0 \leq \beta \leq 2)$ via a sum-of-exponentials approximation efficiently on the interval $[\delta, T]$ with $\delta = \min_{1 \leq K \leq N} \Delta t_n$ and the absolute error $\epsilon$. That is to say, there exit positive real numbers $s_i$ and $\omega_i (i = 1, \ldots, K_{exp})$ such that

$$\left| \frac{1}{t^\beta} - \sum_{i=1}^{K_{exp}} \omega_i e^{-s_i t} \right| \leq \epsilon, t \in [\delta, T].$$

Then we split the derivative into a sum of local part and history part as following:

$$\begin{aligned}
{}_0^C D_t^\alpha \phi(\mathbf{x}, t_{k+1}) &= \frac{1}{\Gamma(1-\alpha)} \int_0^t \frac{\phi_s(\mathbf{x}, s)}{(t-s)^\alpha} ds, \\
&= \frac{1}{\Gamma(1-\alpha)} \int_{t_k}^{t_{k+1}} \frac{\phi_s(\mathbf{x}, s)}{(t-s)^\alpha} ds + \frac{1}{\Gamma(1-\alpha)} \int_0^{t_k} \frac{\phi_s(\mathbf{x}, s)}{(t-s)^\alpha} ds \\
&:= C_l(t_{k+1}) + C_h(t_{k+1}),
\end{aligned}$$

The local part $C_l(t_{k+1})$ and the history part $C_h(t_{k+1})$ are approximated respectively as follows:

$$C_l(t_{k+1}) \approx \frac{\phi(\mathbf{x}, t_{k+1}) - \phi(\mathbf{x}, t_k)}{\Delta t^\alpha \Gamma(2-\alpha)},$$

and

$$\begin{aligned}
C_h(t_{k+1}) &= \frac{1}{\Gamma(1-\alpha)} \left[ \frac{\phi(\mathbf{x}, t_k)}{\Delta t^\alpha} - \frac{\phi(\mathbf{x}, t_0)}{t_{k+1}^\alpha} - \alpha \int_0^{t_k} \frac{\phi(\mathbf{x}, s)}{(t-s)^{1+\alpha}} ds \right] \\
&\approx \frac{1}{\Gamma(1-\alpha)} \left[ \frac{\phi(\mathbf{x}, t_k)}{\Delta t^\alpha} - \frac{\phi(\mathbf{x}, t_0)}{t_{k+1}^\alpha} - \alpha \sum_{i=1}^{K_{exp}} \omega_i \int_0^{t_k} e^{-(t_{k+1}-\tau)s_i} \phi(\mathbf{x}, \tau) d\tau \right] \\
&= \frac{1}{\Gamma(1-\alpha)} \left[ \frac{\phi(\mathbf{x}, t_k)}{\Delta t^\alpha} - \frac{\phi(\mathbf{x}, t_0)}{t_{k+1}^\alpha} - \alpha \sum_{i=1}^{K_{exp}} \omega_i \mathbf{U}_{hist,i}(t_{k+1}) \right],
\end{aligned}$$

with

$$\begin{aligned}
\mathbf{U}_{hist,i}(t_{k+1}) &= \int_0^{t_k} e^{-(t_{k+1}-\tau)s_i}\phi(\mathbf{x},\tau)d\tau \\
&= e^{-s_i\Delta t}\mathbf{U}_{hist,i}(t_k) + \int_{t_{k-1}}^{t_k} e^{-(t_{k+1}-\tau)s_i}\phi(\mathbf{x},\tau)d\tau \\
&\approx e^{-s_i\Delta t}\mathbf{U}_{hist,i}(t_k) + \frac{e^{-s_i\Delta t}}{s_i^2\Delta t}\Big[(e^{-s_i\Delta t} - 1 + s_i\Delta t)\phi(\mathbf{x},t_k) \\
&\quad + (1 - e^{-s_i\Delta t} - e^{-s_i\Delta t}s_i\Delta t)\phi(\mathbf{x},t_{k-1})\Big].
\end{aligned}$$

Denote

$$_0^C\mathbf{D}_t^\alpha \phi^{n+1} = \frac{\phi^{n+1} - \phi^n}{\Delta t^\alpha \Gamma(2-\alpha)} + \frac{1}{\Gamma(1-\alpha)}\Big[\frac{\phi^n}{\Delta t^\alpha} - \frac{\phi^0}{t_{n+1}^\alpha} - \alpha \sum_{i=1}^{K_{exp}} \omega_i \mathbf{U}_{hist,i}^{n+1}\Big],$$

where

$$\begin{aligned}
\mathbf{U}_{hist,i}^{n+1} &= e^{-s_i\Delta t}\mathbf{U}_{hist,i}^n + \frac{e^{-s_i\Delta t}}{s_i^2\Delta t}\Big[(e^{-s_i\Delta t} - 1 + s_i\Delta t)\phi^n \\
&\quad + (1 - e^{-s_i\Delta t} - e^{-s_i\Delta t}s_i\Delta t)\phi^{n-1}\Big].
\end{aligned}$$

Then the time discrete scheme of (2.1) based on the fast algorithm is:

**Scheme 3.1.** *Set the initial step $\phi^0 = \phi(\mathbf{x},0)$, and calculate $\phi^1$ via a first-order scheme.*

$$_0^C\mathbf{D}_t^\alpha \phi^1 = -M\varepsilon^2 \Delta^2 \phi^1 - \nabla \cdot ((|\nabla \phi^0|^2 - 1)\nabla \phi^0). \tag{3.2}$$

*Then, after obtained $\phi^{n-1}$ and $\phi^n$, $\forall n \geq 1$, we can calculate $\phi^{n+1}$ via*

$$_0^C\mathbf{D}_t^\alpha \phi^{n+1} = -M\varepsilon^2 \Delta^2 \phi^{n+1} - \nabla \cdot (2(|\nabla \phi^n|^2 - 1)\nabla \phi^n - (|\nabla \phi^{n-1}|^2 - 1)\nabla \phi^{n-1}). \tag{3.3}$$

Next, we will propose time discretization and spatial discretization to obtain linear and discrete numerical approximations.

### 3.2 Full Discretization using Spectral Method

In this paper, fourier spectral method is employed to handle the spatial discretization, since it is one of the most suitable spatial approximation methods for periodic problems. We use the following Fourier basis functions:

$$S_K = span\{\exp^{-i\mathbf{k}\mathbf{x}}, -K \leq k_1, k_2 \leq K\}, \tag{3.4}$$

where $\mathbf{k} = (k_1, k_2)$, $K$ is a positive integer, $N = 2K + 1$ are the number of the Fourier modes.

Then the weak scheme formulation for model (2.1) can be written as following:

**Scheme 3.2.** *Give the initial condition $\phi_N^0 = \phi_0$, and calculate $\phi_1^N$ through a first-order scheme.*

$$\left(_0^C\mathbf{D}_t^\alpha \phi_N^1, \varphi_N\right) = -M\left(\varepsilon^2 \Delta^2 \phi_N^1 - \nabla \cdot ((|\nabla \phi_N^0|^2 - 1)\nabla \phi_N^0), \varphi_N\right). \tag{3.5}$$

*After calculated $\phi_N^n, \phi_N^{n-1} \in S_K$, $\forall n \geq 1$, find $\phi_N^{n+1} \in S_K$ such that*

$$\left(_0^C\mathbf{D}_t^\alpha \phi_N^{n+1}, \varphi_N\right) = -M\left(\varepsilon^2 \Delta^2 \phi_N^{n+1} - \nabla \cdot (2(|\nabla \phi_N^n|^2 - 1)\nabla \phi_N^n - (|\nabla \phi_N^{n-1}|^2 - 1)\nabla \phi_N^{n-1}), \varphi_N\right), \tag{3.6}$$

where

$$_0^C\mathbf{D}_t^\alpha \phi_N^{n+1} = \frac{\phi_N^{n+1} - \phi_N^n}{\Delta t^\alpha \Gamma(2-\alpha)} + \frac{1}{\Gamma(1-\alpha)}\left[\frac{\phi_N^n}{\Delta t^\alpha} - \frac{\phi_N^0}{t_{n+1}^\alpha} - \alpha \sum_{i=1}^{K_{exp}} \omega_i \mathbf{U}_{hist,i}^{n+1}\right],$$

and

$$\mathbf{U}_{hist,i}^{n+1} = e^{-s_i \Delta t}\mathbf{U}_{hist,i}^n + \frac{e^{-s_i \Delta t}}{s_i^2 \Delta t}\left[(e^{-s_i \Delta t} - 1 + s_i \Delta t)\phi_N^n \right.$$
$$\left. + (1 - e^{-s_i \Delta t} - e^{-s_i \Delta t} s_i \Delta t)\phi_N^{n-1}\right].$$

Denote the coefficient $\mathbf{k}$th mode Fourier coefficient of the unknowns $\phi_N^{n+1}$, and express the approximate solution $\phi_N^{n+1}$ in form of a truncated Fourier expansion:

$$\phi_N^{n+1}(\mathbf{x}) = \sum_{k_1,k_2=-K}^{K} \hat{\phi}_N^{n+1} \exp^{-i\mathbf{k}\mathbf{x}},$$

Then by applying the Fourier transformation to (3.6), we obtain a number of ordinary differential equations for each mode $\mathbf{k}$ in the Fourier space,

$$_0^C\mathbf{D}_t^\alpha \hat{\phi}_N^{n+1} = -M\varepsilon^2 |\mathbf{k}|^4 \hat{\phi}_N^{n+1} - i\mathbf{k}\left\{(2(|\nabla \phi_N^n|^2 - 1)\nabla \phi_N^n - (|\nabla \phi_N^{n-1}|^2 - 1)\nabla \phi_N^{n-1})\right\}_{\mathbf{k}}, \tag{3.7}$$

where $\mathbf{k} = \sqrt{k_1^2 + k_2^2}$ is the magnitude of $\mathbf{k}$ and $\{f\}_{\mathbf{k}}$ represent the $\mathbf{k}$th-mode Fourier coefficient of the function $f$. The Fourier coefficients of the nonlinear term can be calculated by performing the discrete fast Fourier transform(FFT). The total cost to compute $h_N^{n+1}$ from (3.7) is $O(N^2 \log(N))$ at each time step.

It could be easily shown that the full discrete linear scheme is uniquely solvable, and it also preserve the mass conservation, but energy dissipation properties in the full discrete should leave to the readers who is interest in.

# 4 Numerical Results

In this section, we test our proposed full discrete numerical scheme (3.6), and study the time fractional MBE model numerically. In particular, we conduct several numerical simulations of the time fractional MBE model with different time derivative order $\alpha$. Here we chose periodic boundary conditions in the square domain $[0, 2\pi] \times [0, 2\pi]$. Define the roughness measure function $W(t)$ as follows:

$$W(t) = \sqrt{\frac{1}{|\Omega|} \int_\Omega \left(\phi(x,y,t) - \overline{\phi}(x,y,t)\right)^2 d\Omega}, \tag{4.1}$$

where $\overline{\phi}(x,y,t) = \int_\Omega \phi(x,y,t) d\Omega$.

## 4.1 Accuracy test

We begin with testing the time accuracy for the full discrete Scheme 3.2. Consider time-fractional MBE model (2.1) with the initial condition, namely

$$\phi(\mathbf{x}, 0) = 0.1(\sin 3x \sin 2y + \sin 5x \sin 5y). \tag{4.2}$$

The space is discretized by $128 \times 128$ grid points by the Fourier spectral method, and set $\varepsilon = 1$. We use numerical results of scheme (3.6) with $\delta t = 0.00001$ and $N = 128$ as the exact solution since the exact solution for time-fractional MBE growth model is unknown. the numerical errors are computed at $t = 1$. Figure 4.1 shows the $L^2$-errors versus time step $\triangle t$ for the fractional MBE growth model with $\alpha = 0.5$(left) and $\alpha = 0.9$ (right). The expected $2 - \alpha$ order convergence rate in time is obtained.

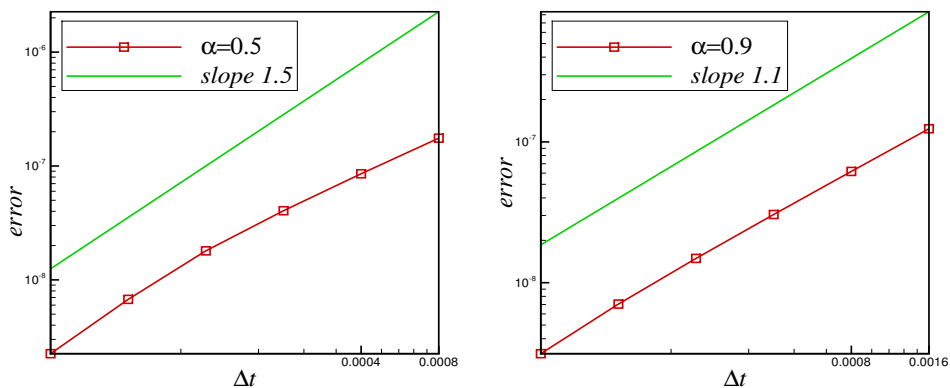

Figure 4.1: Time accuracy convergence test. In the figures, the error versus time step $\triangle t$ for the MBE growth model are shown using scheme (3.6) with $\alpha = 0.5$ (left) and $\alpha = 0.9$ (right), fixed parameter $\varepsilon^2 = 1$.

## 4.2 Example 1

Using the same initial condition as (4.2) and parameters as previous sub-section, we investigate how $\alpha$ affects the evolution dynamics. Here we pick $\alpha$ as 0.3, 0.5, 0.9, 1.0, and use time step $\delta t = 10^{-3}$. The contour lines of numerical solutions of the height function $\phi$ for the time-fractional MBE growth model with different $\alpha$ are shown in Figure 4.2, Figure 4.3, Figure 4.4 and Figure 4.5, respectively. Similar coarsening dynamics are observed, but with different time scales. At $t = 30$, all cases reach the same quasi-steady state. This indicates that the time-fractional MBE model would affects the time scaling dynamics of molecular beam epitaxy, and it doesn't affect its steady state, which is consistent with the idea of fractional time derivative and the findings in other fractional phase field models [26].

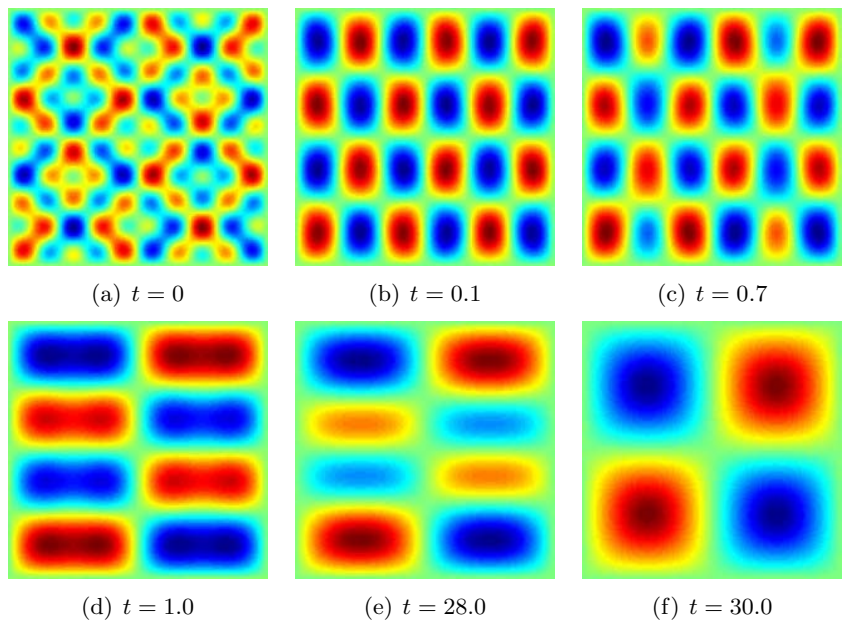

(a) $t = 0$  (b) $t = 0.1$  (c) $t = 0.7$

(d) $t = 1.0$  (e) $t = 28.0$  (f) $t = 30.0$

Figure 4.2: The isolines of numerical solutions of the height function $\phi$ for the MBE growth model using Crank-Nicolson scheme (3.6). With $\varepsilon^2 = 0.1, \alpha = 0.3$ and time step $\delta t = 0.001$. Snapshots are taken at $t = 0, 0.1, 0.7, 1.0, 28, 30$, respectively.

To further resolve how the time scaling is affected by the fractional order $\alpha$, we plot and compare the time evolution of the effective energy and roughness in $t \in [0, 120]$ with different fractional order $\alpha$, as shown in Figure 4.6. In addition, since the dynamic evolves too rapidly between the time interval $[0, 6]$ and $[6, 80]$, we further provide a zoom-view for time evolution of the energy(left) and rough(right) with different $\alpha$ at time region $[0, 6]$(shown in Figure 4.7) and $[6, 80]$(shown Figure 4.8) respectively. As the figures shows that, the energy will be decay faster at beginning as smaller $\alpha$ (shown in Figure 4.7), then it will be decay slower as smaller $\alpha$ (shown in Figure 4.8). And such relation is linearly related. And in the end, they reach the same steady-state energy level.

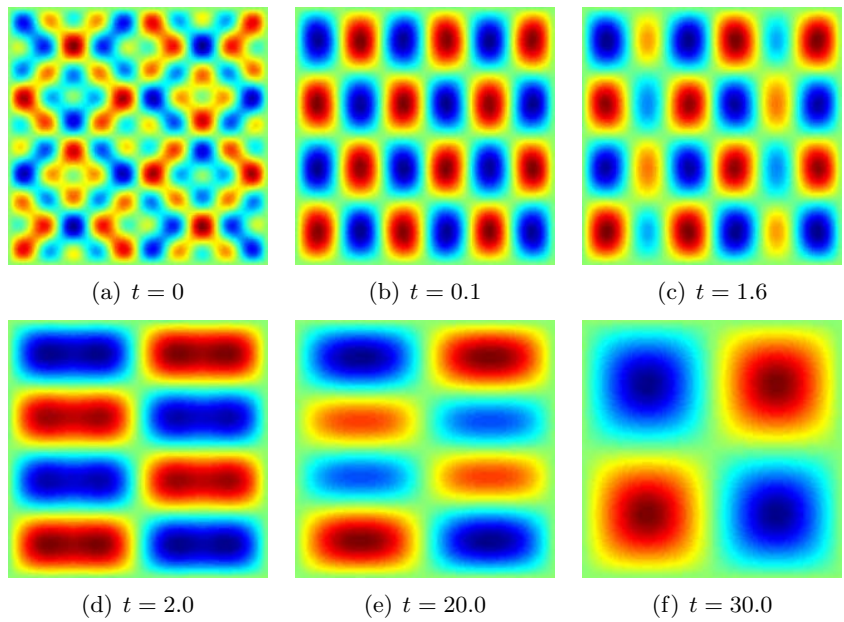

Figure 4.3: The isolines of numerical solutions of the height function $\phi$ for the MBE growth model using Crank-Nicolson scheme (3.6). With $\varepsilon^2 = 0.1, \alpha = 0.5$ and time step $\delta t = 0.001$. Snapshots are taken at $t = 0, 0.1, 1.6, 2.0, 20, 30$, respectively.

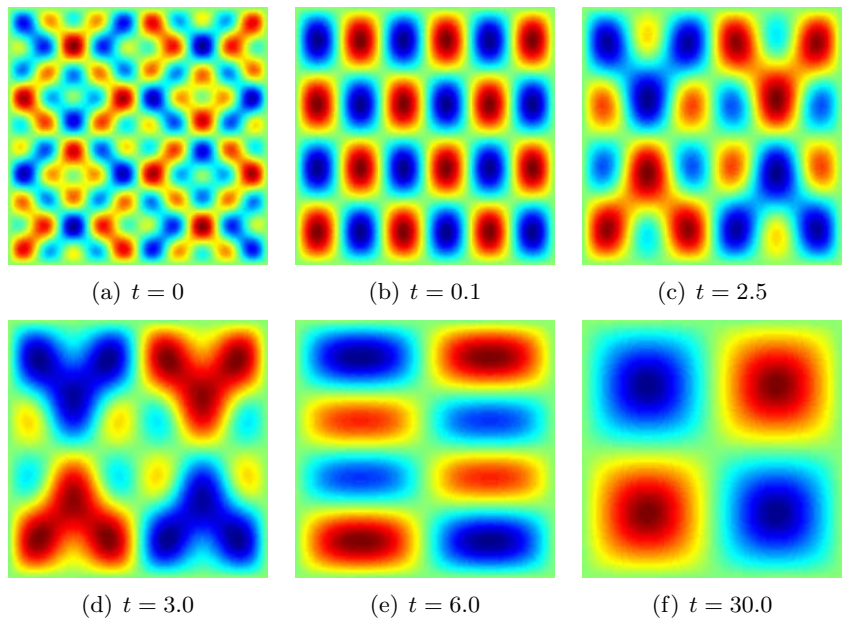

Figure 4.4: The isolines of numerical solutions of the height function $\phi$ for the MBE growth model using Crank-Nicolson scheme (3.6). With $\varepsilon^2 = 0.1, \alpha = 0.9$ and time step $\delta t = 0.001$. Snapshots are taken at $t = 0, 0.1, 2.5, 3, 6, 30$, respectively.

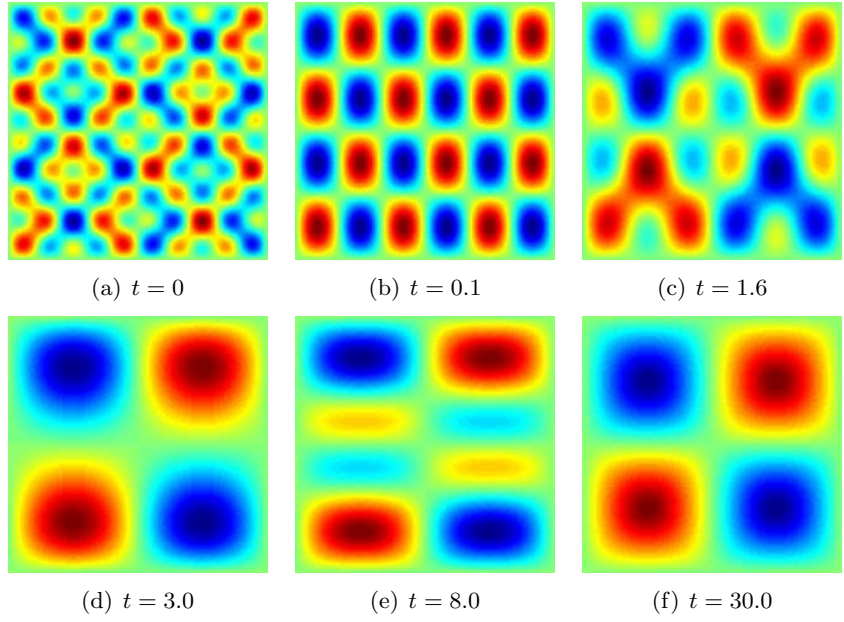

(a) $t=0$    (b) $t=0.1$    (c) $t=1.6$

(d) $t=3.0$    (e) $t=8.0$    (f) $t=30.0$

Figure 4.5: The isolines of numerical solutions of the height function $\phi$ for the MBE growth model using Crank-Nicolson scheme (3.6). With $\varepsilon^2 = 0.1, \alpha = 1.0$ and time step $\delta t = 0.001$. Snapshots are taken at $t = 0, 0.1, 1.6, 3, 8, 30$, respectively.

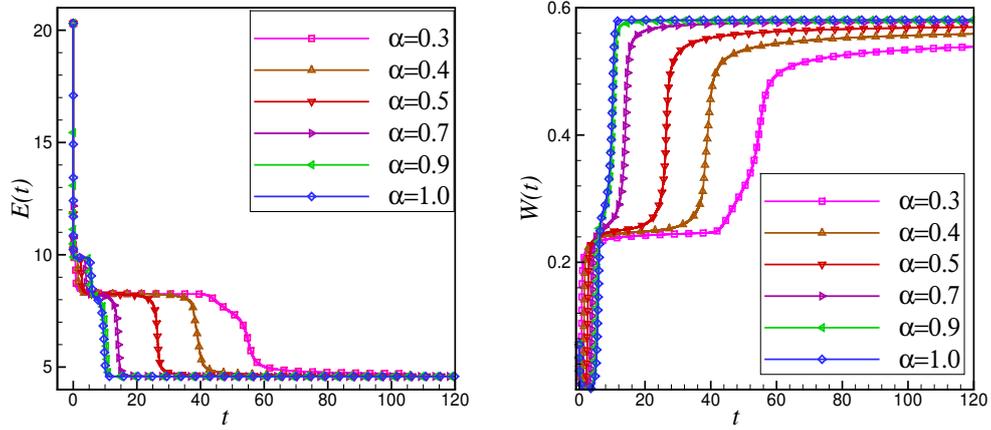

Figure 4.6: Time evolution of the energy(left) and rough(right) for the MBE growth model using scheme (3.6) when $t \in [0, 120]$. The time step is set at $\triangle t = 10^{-3}$ and $\alpha = 0.3, 0.4, 0.5, 0.7, 0.9, 1.0$.

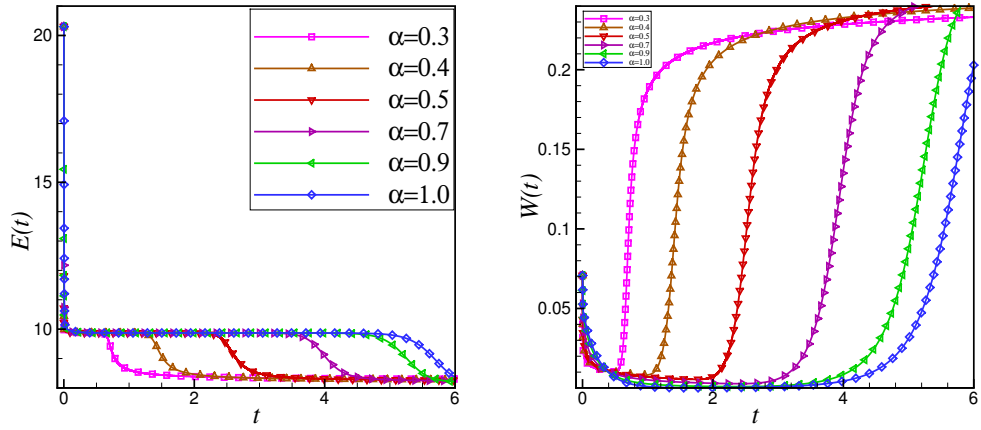

Figure 4.7: Time evolution of the energy(left) and rough(right) for the MBE growth model using scheme (3.6) when $t \in [0, 6]$. The time step is set at $\triangle t = 10^{-3}$ and $\alpha = 0.3, 0.4, 0.5, 0.7, 0.9, 1.0$.

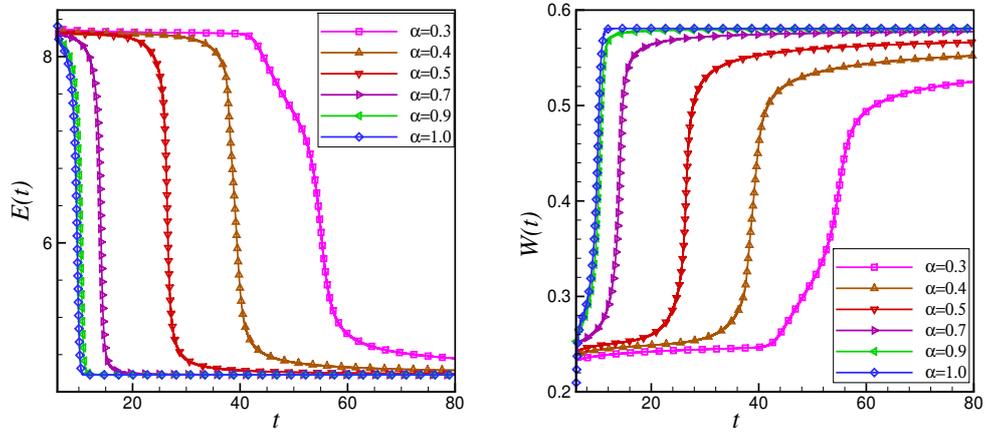

Figure 4.8: Time evolution of the energy(left) and rough(right) for the MBE growth model using scheme (3.6) when $t \in [6, 80]$. The time step is set at $\triangle t = 10^{-3}$ and $\alpha = 0.3, 0.4, 0.5, 0.7, 0.9, 1.0$.

## 4.3 Example II

In this sub-section, we study the coarsening dynamics with random initial condition. It is known for the integer MBE model with slope selection (1.5), the roughness (the standard deviation of the height profile) scales as $O(t^{\frac{1}{3}})$ and the effective energy scales as $O(t^{-\frac{1}{3}})$ [19, 22, 23, 39, 45]. One might wonder how the scaling rate could be manipulated. Here, we will attempt to answer this question numerically.

In this example, we choose the polynomial degree $N = 512$ and the initial condition a random state by assigning a random number which varying from $-0.001$ to $0.001$ to each grid points, then perform numerical simulations of coarsening dynamics. Here the domain is chosen as $[0 \ L_x] \times [0 \ L_y]$ with $L_x = L_y = 2\pi$ and $\varepsilon = 0.03$.

In Figure 4.9, time evolution of the energy (left) and roughness (right) with random number condition when $t \in [0, 500]$ and different $\alpha$ ($\alpha = 0.1, 0.3, 0.5, 0.7, 0.9, 1.0$) are plotted respectively. Also, we plot the time evolution of the energy (left) and roughness (right) when $t \in [0.1, 100]$ with a $\log_{10}$-scale for both x-axis and y-axis in Figure 4.10. By the linear least square fitting, we get the absolute value of slope for each linear line of energy and roughness, $\beta_{E(t)}(\alpha), \beta_{W(t)}(\alpha)$, defined as

$$\log_{10} E(\alpha, t) = \beta_E^0(\alpha) + \beta_E(\alpha) \log_{10} t, \quad \log_{10} W(\alpha, t) = \beta_W^0(\alpha) + \beta_W(\alpha) \log_{10} t. \quad (4.3)$$

They are plotted in Figure 4.11, with the line $\frac{1}{3}\alpha$ plotted for comparison. We observe that the energy decrease approximately as $)(t^{-\frac{1}{3}\alpha})$ and the growth rate of the roughness function is approximately as $O(t^{\frac{1}{3}\alpha})$.

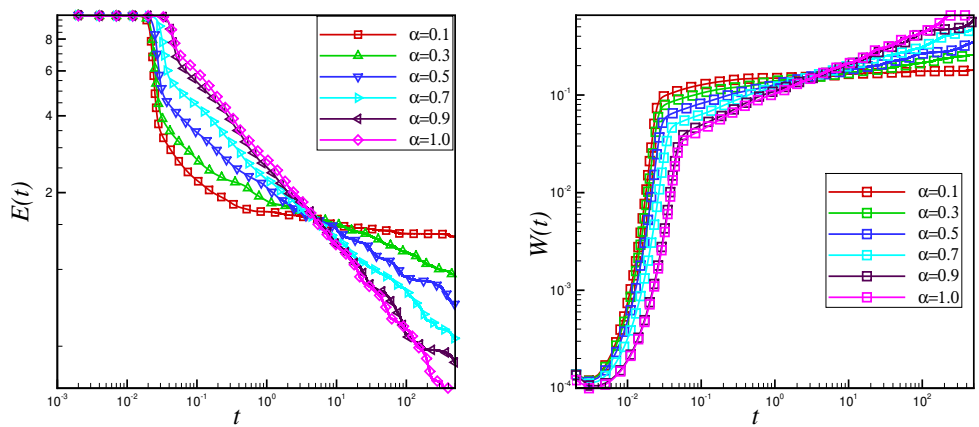

Figure 4.9: Time evolution of the energy (left) and roughness (right) for the MBE growth model using scheme (3.6) with random number condition when $t \in [0.001, 500]$ respectively.

Furthermore, the contour lines of numerical solutions of the height function $\phi$ and its Laplacian $\Delta\phi$ for the time-fractional MBE growth model with random initial condition and $\alpha = 0.1, 0.3, 0.9, 1.0$ are shown in Figure 4.12, 4.13, 4.14 and 4.15, respectively. The time step is $\delta t = 0.001$. Snapshots are taken at $t = 0, 1, 10, 50, 100, 500$. At beginning, the coarsening dynamics will be faster as smaller $\alpha$ but it will be much more slower after reach at time point. This agrees qualitatively well with our results in the first example.

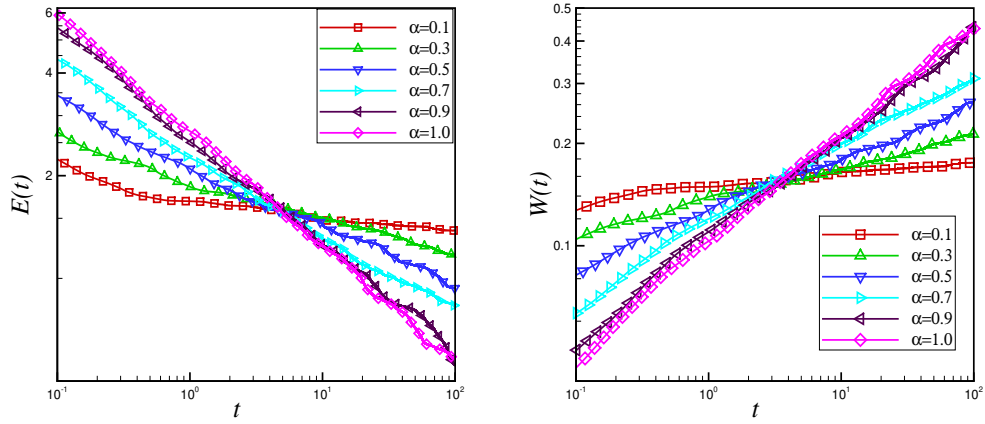

Figure 4.10: Time evolution of the energy (left) and roughness (right) for the MBE growth model using scheme (3.6) with random number condition when $t \in [0.1, 100]$ respectively.

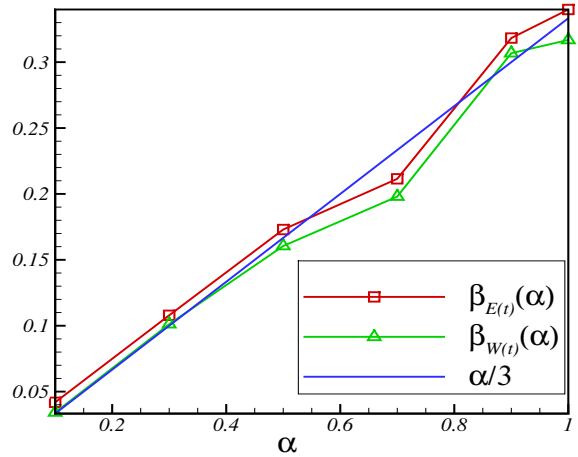

Figure 4.11: Time evolution of the energy (left) and roughness (right) for the MBE growth model using Crank-Nicolson scheme (3.6) with random number condition when $t \in [0.1, 100]$ respectively.

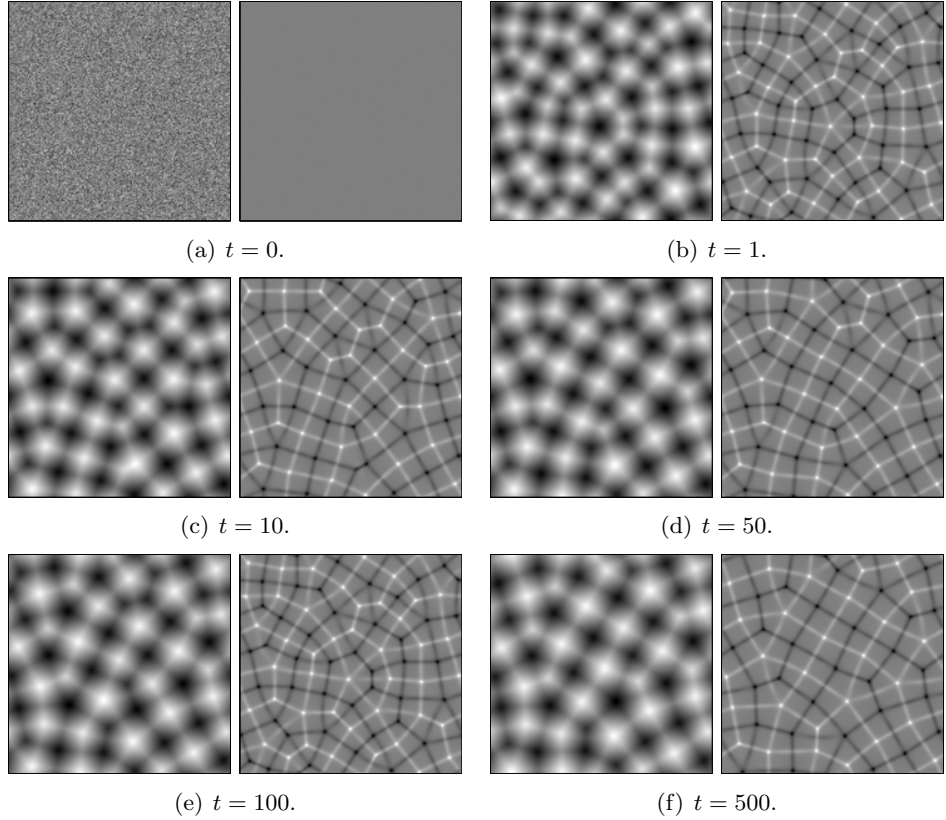

Figure 4.12: The isolines of numerical solutions of the height function $\phi$ and its Laplacian $\Delta \phi$ for the MBE growth model with random initial condition. The time step is $\delta t = 0.001$ and $\alpha = 0.1$. Snapshots are taken at $t = 0, 1, 10, 50, 100, 500$, respectively.

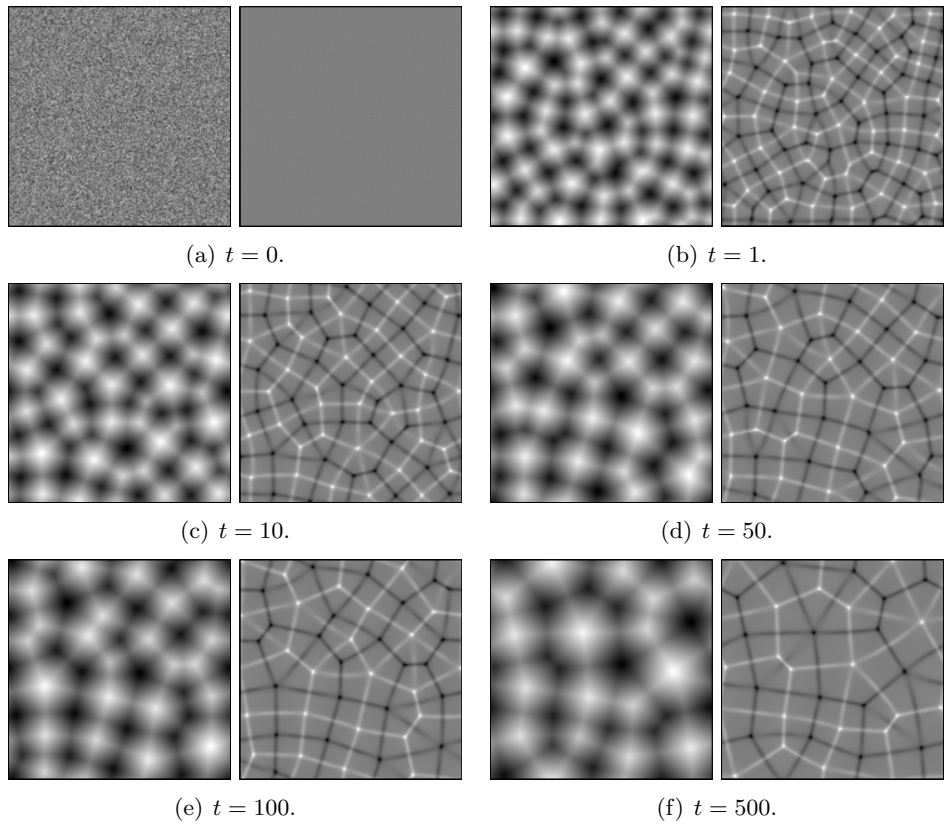

Figure 4.13: The isolines of numerical solutions of the height function $\phi$ and its Laplacian $\Delta\phi$ for the MBE growth model with random initial condition. The time step is $\delta t = 0.001$ and $\alpha = 0.3$. Snapshots are taken at $t = 0, 1, 10, 50, 100, 500$, respectively.

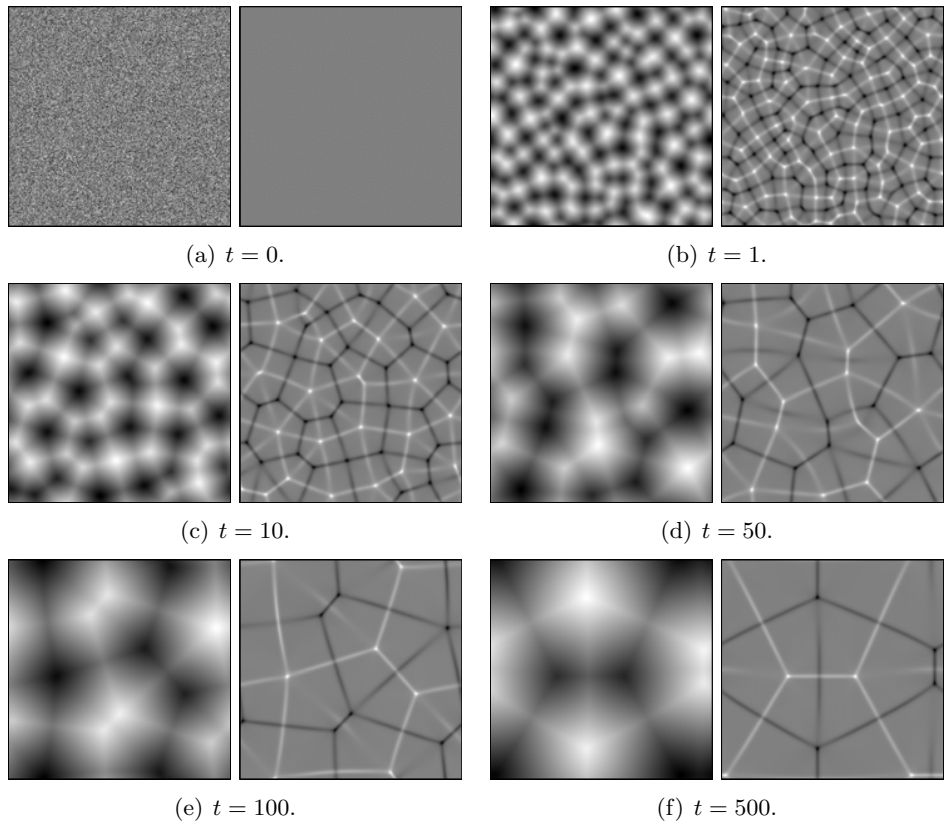

Figure 4.14: The isolines of numerical solutions of the height function $\phi$ and its Laplacian $\Delta\phi$ for the MBE growth model with random initial condition. The time step is $\delta t = 0.001$ and $\alpha = 0.9$. Snapshots are taken at $t = 0, 1, 10, 50, 100, 500$, respectively.

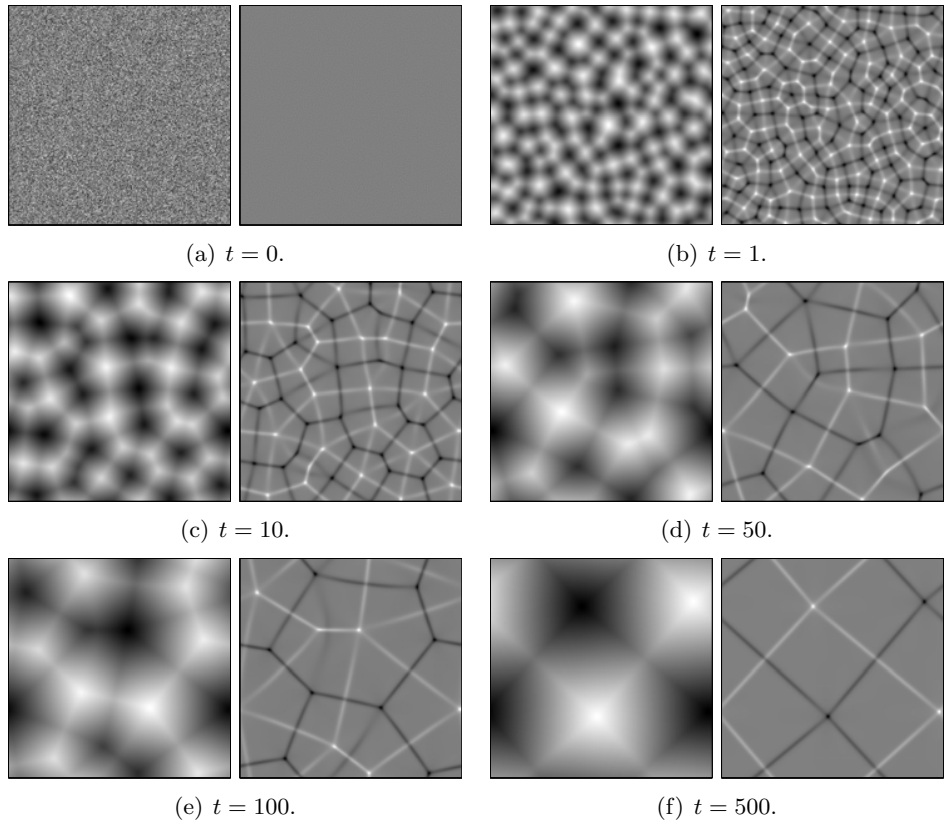

(a) $t = 0$.  (b) $t = 1$.
(c) $t = 10$.  (d) $t = 50$.
(e) $t = 100$.  (f) $t = 500$.

Figure 4.15: The isolines of numerical solutions of the height function $\phi$ and its Laplacian $\Delta\phi$ for the MBE growth model with random initial condition. The time step is $\delta t = 0.001$ and $\alpha = 1.0$. Snapshots are taken at $t = 0, 1, 10, 50, 100, 500$, respectively.

# 5 Conclusion

In this paper, we propose a time fractional molecular beam epitaxy model with slope selection. Then we develop an accurate and efficient numerical scheme for the model. The Caputo fractional derivative was evaluated by the fast algorithm, making the long time dynamic simulation in 2D practical. The space was discretized by the Fourier spectral method, allowing us to capture fine resolution details in the coarsening dynamics of molecular beam epitaxy growth. The mass conservation property is preserved by the full discrete scheme, and the energy dissipation law was checked/verified numerically. We have carried out the convergence accuracy test by the smooth initial condition, showing the $2 - \alpha$ order accuracy in time. In addition, numerical simulation of coarsening dynamics with random initial condition was presented, as well.

Some promising dynamics have been observed for the time-fractional MBE model. In particular, our numerical simulations indicate a scaling law of $O(t^{-\frac{\alpha}{3}})$ for the energy decay and $O(t^{\frac{\alpha}{3}})$ for roughness growth. To our best knowledge, this is the first time report of this observation in literature. As a future project, the theoretical analysis for the energy dissipation properties for the continuous fractional gradient flow models and their numerical schemes will be considered.

## Acknowledgment

Lizhen Chen would like to acknowledge the support from National Science Foundation of China through Grant 11671166 and U1530401 , Postdoctoral Science Foundation of China through Grant 2015M580038. Jia Zhao is partially supported by a faculty startup grant at Utah State University. Waixiang Cao's Research supported in part by NSFC grant No.11501026, and and the Fundamental Research Funds for the Central Universities No. 2017NT10. Hong Wang is partially supported by the OSD/ARO MURI Grant W911NF-15-1-0562 and by the National Science Foundation under Grant DMS-1620194.